\documentclass[10pt]{article}

\usepackage{amssymb,amsfonts, amsopn, latexsym,graphicx,makeidx,enumerate}

\newcommand{\U}{{\mathcal U}}

\newcommand{\C}{{\mathbb C}}
\newcommand{\Z}{{\mathbb Z}}

\newcommand{\D}{{\mathbb D}}

\newcommand{\mf}{{F_{f}}}
\newcommand{\mfo}{{F_{f, \mathbf 0}}}
\newcommand{\mfoo}{{F_{f_0}}}

\newcommand{\dm}{\operatorname{dim}}

\newcommand{\rank}{\operatorname{rank}}
\newcommand{\arrow}[1]{\stackrel{#1}{\longrightarrow}}
\newcommand{\text}[1]{\mbox{\rm {#1}}}

\newcommand{\pc}{{\Gamma^1_{f, z_0}}}

\newtheorem{defn0}{Definition}[section]
\newtheorem{prop0}[defn0]{Proposition}
\newtheorem{conj0}[defn0]{Conjecture}
\newtheorem{thm0}[defn0]{Theorem}
\newtheorem{lem0}[defn0]{Lemma}
\newtheorem{corollary0}[defn0]{Corollary}
\newtheorem{example0}[defn0]{Example}
\newtheorem{remark0}[defn0]{Remark}
\newtheorem{question0}[defn0]{Question}

\newenvironment{defn}{\begin{defn0}\hskip -.06in .}{\end{defn0}}

\newenvironment{thm}{\begin{thm0}\hskip -.06in .}{\end{thm0}}
\newenvironment{lem}{\begin{lem0}\hskip -.06in .}{\end{lem0}}
\newenvironment{cor}{\begin{corollary0}\hskip -.06in .}{\end{corollary0}}
\newenvironment{exm}{\begin{example0}\hskip -.06in .\rm}{\end{example0}}
\newenvironment{rem}{\begin{remark0}\hskip -.06in .\rm}{\end{remark0}}
\newenvironment{ques}{\begin{question0}\hskip -.06in .\rm}{\end{question0}}

\newcommand{\defref}[1]{Definition~\ref{#1}}

\newcommand{\thmref}[1]{Theorem~\ref{#1}}
\newcommand{\lemref}[1]{Lemma~\ref{#1}}
\newcommand{\corref}[1]{Corollary~\ref{#1}}

\newcommand{\quesref}[1]{Question~\ref{#1}}

\newcommand{\qed}{\mbox{$\Box$}}
\newenvironment{proof}{\noindent {\bf Proof.}}{\qed\vskip 6pt}

\title{Hypersurface Singularities and Milnor Equisingularity\footnote{The second author would like to thank the Abdus Salam ICTP for their hospitality; most of this paper was written during a visit there.
\newline   AMS subject classifications 32B15, 32C35, 32C18, 32B10.
\newline  keywords: hypersurface singularity, Milnor fiber, swing, polar curve, vanishing cycles, discriminant, Cerf diagram, intersection diagram, perverse sheaves}}

\author{L\^e D\~ung Tr\'ang and David B. Massey}

\date{}

\begin{document}

\baselineskip= 14pt
\maketitle

\begin{abstract} Suppose that $f$ defines a singular, complex affine hypersurface. If the critical locus of $f$  is one-dimensional at the origin, we obtain new general bounds on the ranks of the homology groups of the Milnor fiber, $F_{f, \mathbf 0}$, of $f$ at the origin, with either integral or $\mathbb Z/p\mathbb Z$ coefficients. If the critical locus of $f$ has arbitrary dimension, we show that the smallest possibly non-zero reduced Betti number of $F_{f, \mathbf 0}$ completely determines if $f$ defines a family of isolated singularities, over a smooth base, with constant Milnor number. This result has a nice interpretation in terms of the structure of the vanishing cycles as an object in the perverse category.
\end{abstract}

\sloppy




\section{Introduction}\label{sec:intro}

Let $\mathcal U$ be an open neighborhood of the origin in $\C^{n+1}$, let $f:(\mathcal U, \mathbf 0)\rightarrow (\C, 0)$ be complex analytic, and let $s:=\dm_{\mathbf 0}\Sigma f$. 

We will use $\mathbf x:=(x_0, \dots, x_n)$ to denote the standard coordinate functions on $\C^{n+1}$. We will use $\mathbf z:=(z_0, \dots, z_n)$ to denote arbitrary analytic local coordinates on $\mathcal U$ near the origin. All of our constructions and results will depend only on the linear part of the coordinates $\mathbf z$; hence, when we say that the $\mathbf z$ are chosen generically, we mean that the linear part of $\mathbf z$ consists of a generic linear combination of $\mathbf x$ (generic in $\operatorname{PGL}(\mathbb C^{n+1})$).

\smallskip

Let $\mf =\mfo$ denote the Milnor fiber of $f$ at the origin. 
It is well-known (see \cite{katomatsu}) that the reduced integral homology, $\widetilde H_k(\mf)$, of $\mf$ can be non-zero only for $n-s\leq k\leq n$, and is free Abelian in degree $n$.  Cohomologically, this means that $\widetilde H^k(\mf)$ can be non-zero only for $n-s\leq k\leq n$, and is free Abelian in degree $n-s$.

For $s>0$ and arbitrary $f$, it is not known how to calculate, algebraically, the groups $\widetilde H_{*}(\mf)$ or their ranks. Even for $s=1$, there is no effective, general method for calculating  the ranks of $\widetilde H_{n-1}(\mf)$ and  $\widetilde H_n(\mf)$. However, there are a number of known bounds on the Betti numbers of $F_f$; we need to describe one of these bounds.

For each $s$-dimensional component, $\nu$, of $\Sigma f$, for a generic point $p\in\nu$, for a generic codimension $s$ (in $\mathcal U$) affine linear subspace, $N$ (a normal slice),  containing $p$, the function $f_{|_N}$ has an isolated critical point at $p$ and the Milnor number at $p$ is independent of the choices;  we let ${\stackrel{\circ}{\mu}}_\nu$ denote this common value. 

If the coordinates $(z_0, \dots, z_{s-1})$ are such that $f_{|_{V(z_0, \dots, z_{s-1})}}$ has an isolated critical point at the origin, then the $s$-dimensional L\^e number \cite{lecycles}, $\lambda^s_{f, \mathbf z}(\mathbf 0)$,  at the origin is defined, and $\lambda^s_{f, \mathbf z}(\mathbf 0)=\sum_\nu {\stackrel{\circ}{\mu}}_\nu\big(\nu\cdot V(z_0, \dots, z_{s-1})\big)_{\mathbf 0}$, where the sum is over the $s$-dimensional components $\nu$ of $\Sigma f$. If the coordinates $(z_0, \dots, z_{s-1})$ are sufficiently generic, then $\lambda^s_{f, \mathbf z}(\mathbf 0)$ obtains its minimum value of $\sum_\nu {\stackrel{\circ}{\mu}}_\nu{\operatorname{mult}}_{\mathbf 0}\nu$; we denote this generic value by $\lambda^s_{f}(\mathbf 0)$ (with no subscript by the coordinates). Theorem 3.3 of \cite{lecycles} implies that $\tilde b_{n-s}:=\operatorname{rank}\widetilde H^{n-s}(F_f)\leq \lambda^s_{f}(\mathbf 0)$.

\vskip .3in

We wish to consider families of singularities. Fix a set of local coordinates $\mathbf z$ for $\mathcal U$ at the origin. Let $G:=(z_0, \dots, z_{s-1})$. If $\mathbf q\in\mathcal U$, we define $f_{\mathbf q}:=f_{|_{G^{-1}(G(\mathbf q))}}$.

\bigskip

\begin{defn}\label{def:simplemuconst}
We say that  $f_{\mathbf q}$  is  a simple $\mu$-constant family  at the origin if and only if, at the origin,   $f_{\mathbf 0}$ has an isolated critical point, $\Sigma f$ is smooth, $G_{|_{\Sigma f}}$ has a regular point and, for all $\mathbf q\in\Sigma f$ close to the origin, the Milnor number $\mu_{\mathbf q}(f_{\mathbf q})$ is independent of $\mathbf q$.
\end{defn}

\bigskip

Our interest in simple $\mu$-constant families stems from the fact that they have many ``equisingularity'' properties; see \thmref{thm:milnorequi}. In particular, if $n-s\neq 2$ and $f_{\mathbf q}$  is  a simple $\mu$-constant family  at the origin, then the main theorem of L\^e and Ramanujam in \cite{leramanujam} implies that the local, ambient, topological-type of $V(f_{\mathbf q})$ at $\mathbf q$ is independent of the point $\mathbf q\in\Sigma f$ near the origin.

\vskip .3in

We can now state our main theorem, which tells us that the rank of $\widetilde H_{n-s}(F_f)$ completely determines whether or not $f$ defines a simple $\mu$-constant family.

\bigskip

\noindent{\bf Main Theorem} (\thmref{thm:main}). {\it Suppose that $\dm_{\mathbf 0}\Sigma(f_{\mathbf 0})=0$. 

Then, the rank, $\tilde b_{n-s}$, of $\widetilde H_{n-s}(F_f)$ is equal to $\lambda^s_{f, \mathbf z}(\mathbf 0)$ if and only if $f_{\mathbf q}$ is a simple $\mu$-constant family.
}

\vskip .4in

This general case of the Main Theorem actually follows quickly from the $1$-dimensional case;

\bigskip

\noindent{\bf Theorem} (\thmref{thm:mainone}). {\it Suppose that $\dm_{\mathbf 0}\Sigma f=1$,  and $\dm_{\mathbf 0}\Sigma (f_{|_{V(z_0)}})=0$. Then, the following are equivalent:

\vskip .1in

\noindent a) $f_{\mathbf q}:= f_{|_{z_0^{-1}(z_0(\mathbf q))}}$ is a simple $\mu$-constant family;

\vskip .1in

\noindent b) $\rank\widetilde H_{n-1}(\mf)=\lambda^1_{f, z_0}(\mathbf 0)$;

\vskip .1in

\noindent c) there exists a prime $p$ such that $\dm \widetilde H_{n-1}(\mf;\ \Z/p\Z)=\lambda^1_{f, z_0}(\mathbf 0)$.

\vskip .1in

Thus, if we are not in the Milnor equisingular case, $\rank\widetilde H_{n-1}(\mf)<\lambda^1_{f, z_0}(\mathbf 0) = \lambda^1_{f}(\mathbf 0)$, and this inequality holds with $\Z/p\Z$ coefficients.}

\bigskip

As a corollary to our Main Theorem, we show that it implies that the vanishing cycles of $f$, as an object in the category of perverse sheaves, cannot be semi-simple in non-trivial cases.

\bigskip

In the final section of this paper, we make some final remarks and present counterexamples to some conceivable ``improvements'' on the statement of the Main Theorem.

\section{Milnor Equisingularity}\label{sec:milequi}

There are other conceivable definitions of what one might wish to call a ``simple'' $\mu$-constant family. The definition that we use in \defref{def:simplemuconst} may seem too strong; we used this strong characterization so that it would be clear in the Main Theorem that the condition $\tilde b_{n-s} = \lambda^s_{f, \mathbf z}(\mathbf 0)$ implies that we are in a very trivial case.

In this section, we will show that all other reasonable concepts of $\mu$-constant families are equivalent. Most, if not all, of the equivalences that we prove here can be found in the literature, at least in special cases.

\bigskip

Suppose that $\dm_{\mathbf 0}\Sigma(f_{\mathbf 0})=0$. Then, the analytic cycle $\displaystyle\left[V\Big(z_0, \dots, z_{s-1}, \frac{\partial f}{\partial z_s}, \dots, \frac{\partial f}{\partial z_n}\Big)\right]$ has the origin as a $0$-dimensional component, and $[\mathbf 0]$ appears in this cycle with multiplicity $\mu_{\mathbf 0}(f_{\mathbf 0})$. Thus, at the origin, $C:=\displaystyle\left[V\Big( \frac{\partial f}{\partial z_s}, \dots, \frac{\partial f}{\partial z_n}\Big)\right]$ is purely $s$-dimensional and is properly intersected by $[V(z_0, \dots, z_{s-1})]$. Let $\Gamma^s_{f, \mathbf z}$ denote the sum of the components of $C$ which are not contained in $\Sigma f$, and let $\Lambda^s_{f, \mathbf z}:= C-\Gamma^s_{f, \mathbf z}$. The cycles $\Gamma^s_{f, \mathbf z}$and $\Lambda^s_{f, \mathbf z}$  are, respectively, the $s$-dimensional polar cycle and $s$-dimensional L\^e cycle; see \cite{lecycles}. It follows at once that
$$
\mu_{\mathbf 0}(f_{\mathbf 0}) = \big(\Gamma^s_{f, \mathbf z}\cdot V(z_0, \dots, z_{s-1})\big)_{\mathbf 0}+\big(\Lambda^s_{f, \mathbf z}\cdot V(z_0, \dots, z_{s-1})\big)_{\mathbf 0}.$$

Note that $\Gamma^s_{f, \mathbf z}=0$ is equivalent to the equality of sets $\displaystyle \Sigma f= V\Big(\frac{\partial f}{\partial z_s}, \dots, \frac{\partial f}{\partial z_n}\Big)$.

Using our notation from the introduction, $\Lambda^s_{f, \mathbf z}=\sum_\nu {\stackrel{\circ}{\mu}}_\nu[\nu]$, where the sum is over the $s$-dimensional components $\nu$ of $\Sigma f$, and, by definition, $\lambda^s_{f, \mathbf z}(\mathbf 0)= \big(\Lambda^s_{f, \mathbf z}\cdot V(z_0, \dots, z_{s-1})\big)_{\mathbf 0}$.  Therefore, we obtain:

\begin{lem}\label{lem:dagger} Suppose that $\dm_{\mathbf 0}\Sigma(f_{\mathbf 0})=0$. Then, 
$$
\mu_{\mathbf 0}(f_{\mathbf 0}) =   \big(\Gamma^s_{f, \mathbf z}\cdot V(z_0, \dots, z_{s-1})\big)_{\mathbf 0}+ \sum_\nu {\stackrel{\circ}{\mu}}_\nu\big(\nu\cdot V(z_0, \dots, z_{s-1})\big)_{\mathbf 0}=   \big(\Gamma^s_{f, \mathbf z}\cdot V(z_0, \dots, z_{s-1})\big)_{\mathbf 0}+ \lambda^s_{f, \mathbf z}(\mathbf 0),
$$
where the sum is over all $s$-dimensional components, $\nu$, of $\Sigma f$.

In particular, $\mu_{\mathbf 0}(f_{\mathbf 0}) =  \lambda^s_{f, \mathbf z}(\mathbf 0)$ if and only if $\Gamma^s_{f, \mathbf z}=0$.
\end{lem}

Note that, while $\lambda^s_{f, \mathbf z}(\mathbf 0)$ is not independent of the choice of $\mathbf z$, $\Lambda^s_{f, \mathbf z}$ is, and this fact is very useful. Let  $(\hat z_0, \dots, \hat z_n)$ be a set of local analytic coordinates for $\mathcal U$ which are close to the coordinates $\mathbf z$; let $\hat f_{\mathbf q}$ denote the corresponding analytic family. As  $\dm_{\mathbf 0}\Sigma(f_{\mathbf 0})=0$, $\dm_{\mathbf 0}\Sigma(\hat f_{\mathbf 0})=0$. Let $\widehat C:=\displaystyle\left[V\Big( \frac{\partial f}{\partial \hat z_s}, \dots, \frac{\partial f}{\partial \hat z_n}\Big)\right]$.  Then, Proposition 8.2.a of \cite{fulton} implies that $\lim_{\hat\mathbf z\rightarrow\mathbf z} \widehat C\leq C$, i.e., $\lim_{\hat\mathbf z\rightarrow\mathbf z} \big(\Gamma^s_{f, \hat\mathbf z}+\Lambda^s_{f, \hat\mathbf z}\big)\leq \big(\Gamma^s_{f, \mathbf z}+\Lambda^s_{f, \mathbf z}\big)$. As $\Lambda^s_{f, \mathbf z}$ is independent of the coordinates, we conclude that $\lim_{\hat\mathbf z\rightarrow\mathbf z}\Gamma^s_{f, \hat\mathbf z}\leq \Gamma^s_{f, \mathbf z}$. 

\smallskip

It follows immediately that

\smallskip

\begin{lem}\label{lem:denselem} If there exist coordinates $\mathbf z$ such that $\dm_{\mathbf 0}\Sigma(f_{\mathbf 0})=0$ and $\Gamma^s_{f, \mathbf z}=0$, then the set of coordinates $\hat\mathbf z$ such that $\dm_{\mathbf 0}\Sigma(\hat f_{\mathbf 0})=0$ and $\Gamma^s_{f, \hat\mathbf z}=0$ form an open dense set (again, we mean that the linear portion of $\hat\mathbf z$ is obtained from the standard coordinates by applying a transformation from an open dense set of $\operatorname{PGL}(\mathbb C^{n+1})$).
\end{lem}

\bigskip

There is one more piece of preliminary notation that we need. Consider the blow-up of $\mathcal U$ along the Jacobian ideal, $J(f)$ of $f$, i.e., $B:={\operatorname{Bl}}_{J(f)}\mathcal U$. This blow-up naturally sits inside $\mathcal U\times \mathbb P^n$. Thus, the exceptional divisor $E$ of the blow-up is a cycle in $\mathcal U\times \mathbb P^n$. 

\bigskip

We now give {\bf many} equivalent characterizations of $\mu$-constant families. 

\bigskip

\begin{thm}\label{thm:milnorequi} Let $\mathbf z$ be local coordinates for $\mathcal U$ at the origin such that $\dm_{\mathbf 0}\Sigma(f_{\mathbf 0})=0$. Then, the following are equivalent: \begin{enumerate}[1{.}]
\item For all $\mathbf q\in\Sigma f$ near the origin, $\mu_{\mathbf 0}(f_{\mathbf 0}) = \mu_{\mathbf q}(f_{\mathbf q})$.
\item $\mu_{\mathbf 0}(f_{\mathbf 0})=\lambda^s_f(\mathbf 0)$.
\item $f_{\mathbf q}$ is a simple $\mu$-constant family.
\item $\mu_{\mathbf 0}(f_{\mathbf 0})=\lambda^s_{f, \mathbf z}(\mathbf 0)$.
\item $\Gamma^s_{f, \mathbf z}=0$.
\end{enumerate}

Futhermore, if $n-s\neq 2$, then 1), 2), 3), 4), and 5) above hold if and only if the local, ambient, topological-type of $V(f_{\mathbf q})$ at $\mathbf q$ is independent of the point $\mathbf q\in\Sigma f$ near the origin.

\vskip .3in

\noindent In addition, the following are equivalent:
\begin{enumerate}[a{.}]
\item There exist coordinates $\mathbf z$ such that 1), 2), 3), 4), and 5) above hold.

\item Near the origin, $\Sigma f$ is smooth and $(\mathcal U-\Sigma f, \Sigma f)$ is an $a_f$ stratification, i.e., for all $\mathbf p\in\Sigma f$ near the origin,  for every limiting tangent space, $\mathcal T_{\mathbf p}$, from level hypersurfaces of $f$ approaching $\mathbf p$ , $T_{\mathbf p}(\Sigma f)\subseteq \mathcal T_{\mathbf p}$.

\item $\Sigma f$ is smooth at the origin, and over an open neighborhood of the origin, the exceptional divisor, $E$, as a set, is equal to the projectivized conormal variety to $\Sigma f$ and, hence, as cycles $E=\mu\left[T^*_{\Sigma f}\mathcal U\right]$.

\item For generic $\hat\mathbf z$, $\Gamma^s_{f, \hat\mathbf z} =0$ near the origin.

\item $\Sigma f$ is smooth at the origin and, for all local coordinates $\hat\mathbf z$ such that $V(\hat z_0, \dots, \hat z_{s-1})$ transversely intersects $\Sigma f$ at the origin, $\hat f_{\mathbf q}$ is a simple $\mu$-constant family.

\end{enumerate}

\end{thm}
\begin{proof} Throughout, we work in a sufficiently small neighborhood of the origin. The theorem is stupidly true if $s=0$; so, we suppose that $s\geq 1$.

\medskip

Suppose that 1) holds. Then, \lemref{lem:dagger} implies that 
$$\mu_{\mathbf 0}(f_{\mathbf 0}) =  \big(\Gamma^s_{f, \mathbf z}\cdot V(z_0, \dots, z_{s-1})\big)_{\mathbf 0}+ \sum_\nu \mu_{\mathbf 0}(f_{\mathbf 0})\big(\nu\cdot V(z_0, \dots, z_{s-1})\big)_{\mathbf 0}.$$
Thus, $\Gamma^s_{f, \mathbf z}=0$, and $\sum_\nu \big(\nu\cdot V(z_0, \dots, z_{s-1})\big)_{\mathbf 0}$ must equal $1$, i.e., $\Sigma f$ must a have single smooth $s$-dimensional component, which is transversely intersected by $V(z_0, \dots, z_{s-1})$. Hence, $\lambda^s_f(\mathbf 0)= \sum_\nu {\stackrel{\circ}{\mu}}_\nu\big(\nu\cdot V(z_0, \dots, z_{s-1})\big)_{\mathbf 0}$, and \lemref{lem:dagger} now implies 2).

\medskip

Suppose that 2) holds. Then, \lemref{lem:dagger} implies that $\Gamma^s_{f, \mathbf z}=0$ and, for each $s$-dimensional component $\nu$ of $\Sigma f$, $\big(\nu\cdot V(z_0, \dots, z_{s-1})\big)_{\mathbf 0} = {\operatorname{mult}}_{\mathbf 0}\nu$. In order to conclude 3), we have only to show that $\Sigma f$ must be smooth. 

As $\Gamma^s_{f, \mathbf z}=0$, there is an equality of sets $\Sigma f =V\Big( \frac{\partial f}{\partial z_s}, \dots, \frac{\partial f}{\partial z_n}\Big)$, and so every component of $\Sigma f$ must be at least $s$-dimensional. We conclude that $\Sigma f$ is purely $s$-dimensional. Now, let $\hat\mathbf z$ be a generic choice of coordinates, close to $\mathbf z$. As the Milnor number is upper-semicontinuous, $\mu_{\mathbf 0}(f_{\mathbf 0})\geq \mu_{\mathbf 0}(\hat f_{\mathbf 0})\geq \lambda^s_f(\mathbf 0)$. As $\mu_{\mathbf 0}(f_{\mathbf 0})= \lambda^s_f(\mathbf 0)$, we conclude that  $\mu_{\mathbf 0}(\hat f_{\mathbf 0})= \lambda^s_f(\mathbf 0)$. As prepolar coordinates are generic (see \cite{lecycles}), we may assume that $\hat\mathbf z$ is prepolar. 

	Consider now $g:= f_{|_{V(\hat z_0, \dots, \hat z_{s-2})}}$ (where we mean that $g:=f$ if $s=1$). Then, as $\hat\mathbf z$ is generic, $\Sigma g = \Sigma f\cap V(\hat z_0, \dots, \hat z_{s-2})$ is $1$-dimensional. By induction, Proposition 1.21 of \cite{lecycles} implies that the polar curve $\Gamma^1_{g, z_{s-1}} = \Gamma^s_{f, \mathbf z}\cdot V(\hat z_0, \dots, \hat z_{s-1}) = 0$. Now, Proposition 1.30 of \cite{lecycles} (which uses the non-splitting result, proved independently by Gabrielov \cite{gabrielov}, Lazzeri \cite{lazzeri}, and L\^e \cite{leacampo}) implies that $\Sigma g$ is smooth. As $\hat \mathbf z$ was generic and $\Sigma f$ was purely $s$-dimensional, we conclude that $\Sigma f$ is smooth, and so 3) holds.
	
\medskip

Certainly, 3) implies 1). Therefore, we have shown that 1), 2), and 3) are equivalent. 

\medskip

That 3) implies 4) is immediate. If 4) holds, then easy generalizations of any of the non-splitting arguments of  Gabrielov \cite{gabrielov}, Lazzeri \cite{lazzeri}, and L\^e \cite{leacampo} immediately imply that, at the origin, $\Sigma f$ has a single smooth component and $V(z_0, \dots, z_{s-1})$ transversely intersects $\Sigma f$. Thus, 4) implies 3).

\medskip

Hence, 1), 2), 3) and 4) are equivalent. By \lemref{lem:dagger}, 4) and 5) are equivalent.

\medskip
	
If 3) holds and $n-s\neq 2$, then the classic result of L\^e and Ramanujam in \cite{leramanujam} implies that the local, ambient, topological-type of $V(f_{\mathbf q})$ at $\mathbf q$ is independent of the point $\mathbf q\in\Sigma f$ near the origin.

With no constraint on $n-s$, if the local, ambient, topological-type of $V(f_{\mathbf q})$ at $\mathbf q$ is independent of the point $\mathbf q\in\Sigma f$ near the origin, then 1) holds, since the Milnor number is an invariant of this topological-type.

\medskip

We need to show that a) through e) are equivalent.

\medskip

Assume that a) holds. Then, 3) implies b) by Theorem 6.8 of \cite{lecycles} (which uses the result of L\^e and Saito from \cite{lesaito}).

The equivalence of b) and c) is immediate, and they clearly imply d). That e) implies a) is also clear. It remains for us to show that d) implies e).

Assume d). By \lemref{lem:dagger},  a) holds and, thus, so does c). Hence, $\Sigma f$ is smooth at the origin. Let $\hat\mathbf z$ be such that $V(\hat z_0, \dots, \hat z_{s-1})$ transversely intersects $\Sigma f$ at the origin. Then, c) tells us at once that $\dm_{\mathbf 0}\Sigma(f_{\mathbf 0})=0$ and $\Gamma^s_{f, \hat\mathbf z}=0$, and so, by \lemref{lem:dagger}, a) holds.
\end{proof}

\bigskip

\begin{defn}\label{def:milnorequi}
Whenever the equivalent conditions a), b), c), d), and e) of \thmref{thm:milnorequi} hold, we say that $f$ is Milnor equisingular at the origin.
\end{defn}

\bigskip 

The Main Theorem of this paper, \thmref{thm:main}, tells us that there is another important topological equivalent characterization of Milnor equisingularity. First, however, we must recall some known results and prove the Swing Lemma.

\section{Known Results}\label{sec:known}

We assume that the first coordinate $z_0$ on $\U$ is a generic linear form; in the terminology of \cite{lecycles}, we need for $z_0$ to be ``prepolar'' (with respect to $f$ at the origin). This implies that $\dm_{\mathbf 0}\Sigma(f_{|_{V(z_0)}})\leq s-1$ (provided that $s\neq 0$), that the polar curve, $\pc$, is purely $1$-dimensional at the origin (which vacuously includes the case $\pc=\emptyset$), and $\pc$ has no components contained in $V(f)$ (this last property is immediate in some definitions of the relative polar curve).

\bigskip

For convenience, we assume throughout the remainder of this paper that the neighborhood $\U$ is re-chosen, if necessary, so small that $\Sigma f\subseteq V(f)$, and every component of $\Sigma f$ and $\pc$ contains the origin.

\bigskip

Now, there is the attaching result of L\^e from \cite{leattach} (see, also, \cite{lecycles}): 

\bigskip

\begin{thm}\label{thm:leattach} Up to diffeomorphism, $\mf$ is obtained from $\stackrel{\circ}{\D}\times\mfoo$ by attaching $\tau:=\big(\pc\cdot V(f))_{\mathbf 0}$ handles of index $n$.
\end{thm}

\bigskip

\begin{rem}\label{rem:leattach} On the level of homology, L\^e's attaching result is a type of Lefschetz hyperplane result; it says that, for all $i<n-1$, the inclusion map $\mfoo=\mf\cap V(z_0)\hookrightarrow\mf$ induces isomorphisms $\widetilde H_i(\mfoo)\cong \widetilde H_i(\mf)$, and $\widetilde H_n(\mf)$ and $\widetilde H_{n-1}(\mf)$ are, respectively, isomorphic to the kernel and cokernel of the boundary map $$
\Z^{\tau}\cong H_n(\mf, \mfoo)\arrow{\partial} \widetilde H_{n-1}(\mfoo).
$$

We remind the reader here of the well-known result, first proved by Teissier in \cite{teissiercargese} (in the case of an isolated singularity, but the proof works in general), that
$$
\tau \ =\  \big(\pc\cdot V(f))_{\mathbf 0}\  = \  \Big(\pc\cdot V\left(\frac{\partial f}{\partial z_0}\right)\Big)_{\mathbf 0}+\big(\pc\cdot V(z_0))_{\mathbf 0}.$$
As defined in \cite{lecycles}, the first summand on the right above is $\lambda^0_{f, z_0}(\mathbf 0)$, the $0$-dimensional L\^e number, and second summand on the right above is $\gamma^1_{f, z_0}(\mathbf 0)$, the $1$-dimensional polar number.

\smallskip

If $s=1$, then $\widetilde H_{n-1}(\mfoo)\cong\Z^{\mu_{\mathbf 0}(f_0)}$. Therefore, in the $s=1$ case, one can certainly calculate the difference of the reduced Betti numbers of $\mf$: 
$$
\tilde b_n(\mf)-\tilde b_{n-1}(\mf) = \tau-\mu_{\mathbf 0}(f_0).
$$
Hence, a bound on one of $\tilde b_n(\mf)$ and $\tilde b_{n-1}(\mf)$ automatically produces a bound on the other. As a final comment, it is well-known, and easy to show that $\mu_{\mathbf 0}(f_0)= \gamma_{f, z_0}^1(\mathbf 0)+\lambda_{f, z_0}^1(\mathbf 0)$.

\end{rem}

\vskip .3in

In Proposition 3.1 of \cite{lecycles}, the second author showed how the technique of ``tilting in the Cerf diagram'' or ``the swing'', as used by L\^e and Perron in \cite{leperron} could help refine the result of \thmref{thm:leattach}. Here, we state only the homological implication of Proposition 3.1 of \cite{lecycles}.

\bigskip

\begin{thm}\label{thm:swing}  The boundary map $H_n(\mf, \mfoo)\arrow{\partial} \widetilde H_{n-1}(\mfoo)$ maps a direct summand of $H_n(\mf, \mfoo)$ of rank $\gamma^1$ isomorphically onto a direct summand of $\widetilde H_{n-1}(\mfoo)$.

Thus, the rank of $\widetilde H_n(\mf)$ is at most $\lambda^0_{f, z_0}(\mathbf 0)$, and the rank of $\widetilde H_{n-1}(\mf)$ is at most $\lambda^1_{f, z_0}(\mathbf 0)$.
\end{thm}

\bigskip

However, if one of the components $\nu$ of $\Sigma f$ is itself singular, then the above bounds on the ranks are known not to be optimal. A result of Siersma in \cite{siersmavarlad}, or an easy exercise using perverse sheaves (see the remark at the end of \cite{siersmavarlad}), yields:

\begin{thm}\label{thm:siersmabound} The rank of $\widetilde H_{n-1}(\mf)$ is at most $\sum_\nu {\stackrel{\circ}{\mu}}_\nu$.
\end{thm}

\bigskip

In light of \thmref{thm:swing} and \thmref{thm:siersmabound}, the question is: Is it possible that $\rank\widetilde H_{n-1}(\mf)=\lambda^1_{f, z_0}(\mathbf 0)$? Of course, the answer to this question is ``yes''; if $f$ is Milnor equisingular and $z_0$ is generic, then \thmref{thm:leattach} tells us immediately that $\rank\widetilde H_{n-1}(\mf)=\lambda^1_{f, z_0}(\mathbf 0)$. \thmref{thm:mainone} tells us that the {\bf only} way for $\rank\widetilde H_{n-1}(\mf)$ to equal $\lambda^1_{f, z_0}(\mathbf 0)$ is for $f$ to be Milnor equisingular.

\bigskip

\begin{rem}\label{rem:siersmarem} Of course, if all of the components $\nu$ of $\Sigma f$ are smooth, and $z_0$ is generic, then $\lambda^1_{f, z_0}(\mathbf 0)=\lambda^1_f(\mathbf 0)=\sum_\nu {\stackrel{\circ}{\mu}}_\nu$, and the bounds on the ranks obtained from \thmref{thm:swing} and \thmref{thm:siersmabound} are the same. In addition,  \thmref{thm:siersmabound} is true with arbitrary field coefficients; this yields bounds on the possible torsion in $\widetilde H_{n-1}(\mf)$. We should also remark that the result of Siersma from \cite{siersmavarlad} that we cite above can actually yield a much stronger bound if one knows certain extra topological data -- specifically, one needs that the ``vertical monodromies'' are non-trivial.

Hence, if one of the components of $\Sigma f$ is itself singular , then  $\rank\widetilde H_{n-1}(\mf)<\lambda^1_{f, z_0}(\mathbf 0)$ by  \thmref{thm:siersmabound}. Even in the case where all of the components of $\Sigma f$ are smooth, we could conclude that $\rank\widetilde H_{n-1}(\mf)<\lambda^1_{f, z_0}(\mathbf 0)$ from \cite{siersmavarlad} {\bf if} we knew that one of the vertical monodromies were non-trivial. However, the vertical monodromies are complicated topological data to calculate. In addition, the vertical monodromies can be trivial even when the polar curve is non-empty, i.e., when $f$ is {\bf not} Milnor equisingular.

In \cite{siersmaisoline}, Siersma proved another closely related result. On the level of homology, what he proved was that, if $f$ is not Milnor equisingular, and $\Sigma f$ has a single smooth component, $\nu$, such that ${\stackrel{\circ}{\mu}}_\nu =1$, then $\widetilde H_{n-1}(\mf)=0$; \thmref{thm:mainone}, including the modulo $p$ statement, is a strict generalization of this.

In addition, we should point out that, in \cite{dejong}, Th. de Jong provides evidence that a result like \thmref{thm:mainone} might be true. 
\end{rem}

\bigskip

Before we can prove our Main Theorem, we still need to prove the Swing Lemma.

\section{The Swing}\label{sec:swing}

We prove the one-dimensional version of our Main Theorem by combining the swing technique of \thmref{thm:swing} and the connectivity of the vanishing cycle intersection diagram for isolated singularities, as was proved independently by Gabrielov in \cite{gabrielov} and Lazzeri in \cite{lazzeri}. In some recent notes, M. Tib\u ar uses similar techniques and reaches a number of conclusions closely related to our result.

\bigskip

In Section 3, we referred to the swing (or, tilting in the Cerf diagram), which was used by L\^e and Perron in \cite{leperron} and in Proposition 3.1 of \cite{lecycles}, where the swing was used to prove \thmref{thm:swing}. The swing has also been studied in \cite{caubelthesis}, \cite{tibarbouq}, \cite{lecycles}, \cite{vannierthesis}. As the swing is so crucial to the proof of the Main Theorem, we wish to give a careful explanation of its construction.

\bigskip

Suppose that $\mathcal W$ is an open neighborhood of the origin in $\mathbb C^2$. We will use the coordinates $x$ and $y$ on $\mathcal W$. For notational ease, when we restrict $x$ and $y$ to various subspaces where the domain is clear, we shall continue to write simply $x$ and $y$.

 Let $C$ be a complex analytic curve in $\mathcal W$ such that every component of $C$ contains the origin. We assume that the origin is an isolated point in $V(x)\cap C$ and in $V(y)\cap C$, i.e., that $C$ does not have a component along the $x$- or $y$-axis.
 
 Below, we let $\mathbb D_\epsilon$ denote a closed disk, of radius $\epsilon$, centered at the origin, in the complex plane. We denote the interior of $\mathbb D_\epsilon$ by ${\stackrel{\circ}{\mathbb D}}_{\epsilon}$, and when we delete the origin, we shall superscript with an asterisk, i.e., $\mathbb D^*_\epsilon:= \mathbb D_\epsilon-\{0\}$ and 
 ${\stackrel{{}\hskip -.06in \circ}{\mathbb D^*_{\epsilon}}}:={\stackrel{\circ}{\mathbb D}}_{\epsilon}-\{0\}$.

\bigskip

We select $0<\epsilon_2\ll\epsilon_1\ll 1$ so that:

\medskip

\noindent i): the ``half-open'' polydisk $\mathbb D_{\epsilon_1}\times{\stackrel{\circ}{\mathbb D}}_{\epsilon_2}$ is contained in $\mathcal W$;

 \medskip
 
 \noindent ii):   $(\partial \mathbb D_{\epsilon_1}\times{\stackrel{\circ}{\mathbb D}}_{\epsilon_2})\cap C=\emptyset$ (this uses that the origin is an isolated point in $V(y)\cap C$) ;
 
 \medskip
 
 Note that ii) implies that $(\mathbb D_{\epsilon_1}\times{\stackrel{\circ}{\mathbb D}}_{\epsilon_2})\cap C = ({\stackrel{\circ}{\mathbb D}}_{\epsilon_1}\times{\stackrel{\circ}{\mathbb D}}_{\epsilon_2})\cap C$.
 
 \medskip
 
 \noindent iii):  $\mathbb D_{\epsilon_1}\times {\stackrel{{}\hskip -.06in \circ}{\mathbb D^*_{\epsilon_2}}}\arrow{y} {\stackrel{{}\hskip -.06in \circ}{\mathbb D^*_{\epsilon_2}}}$ is a proper stratified submersion, where the Whitney strata are $\partial \mathbb D_{\epsilon_1}\times {\stackrel{{}\hskip -.06in \circ}{\mathbb D^*_{\epsilon_2}}}$, \hbox{$({\stackrel{\circ}{\mathbb D}}_{\epsilon_1}\times {\stackrel{{}\hskip -.06in \circ}{\mathbb D^*_{\epsilon_2}}})-C$}, and $({\stackrel{\circ}{\mathbb D}}_{\epsilon_1}\times {\stackrel{{}\hskip -.06in \circ}{\mathbb D^*_{\epsilon_2}}})\cap C$.
 
 \medskip
 
 \noindent iv): $({\stackrel{\circ}{\mathbb D}}_{\epsilon_1}\times {\stackrel{{}\hskip -.06in \circ}{\mathbb D^*_{\epsilon_2}}})\cap C\arrow{y} {\stackrel{{}\hskip -.06in \circ}{\mathbb D^*_{\epsilon_2}}}$ is an $m$-fold covering map, where $m:=(C\cdot V(y))_\mathbf 0$.

\bigskip

Let $D:= (\mathbb D_{\epsilon_1}\times{\stackrel{\circ}{\mathbb D}}_{\epsilon_2})\cap (C\cup V(y))$. Let $(x_0, y_0)\in (\mathbb D^*_{\epsilon_1}\times {\stackrel{{}\hskip -.06in \circ}{\mathbb D^*_{\epsilon_2}}})-D$. Let $\sigma:[0,1]\rightarrow\{x_0\}\times {\stackrel{\circ}{\mathbb D}}_{\epsilon_2}$ be a smooth, simple path such that $\sigma(0)= (x_0, y_0)$, $\sigma(1)=:(x_0, y_1)\in C$, and $\sigma([0,1))\subseteq (\{x_0\}\times{\stackrel{\circ}{\mathbb D}}_{\epsilon_2})-D$. 

Let $S$ be the image of $\sigma$; as $\sigma$ is simple, $S$ is homeomorphic to $[0, 1]$. Let $\sigma_0:=y\circ\sigma$ and let $S_0$ be the image of $\sigma_0$. Thus, $S_0$ is homeomorphic to $[0,1]$ and is contained in ${\stackrel{{}\hskip -.06in \circ}{\mathbb D^*_{\epsilon_2}}}$.

\bigskip

\begin{lem}\label{lem:swing}{\rm({\bf The Swing})}  There exists a continuous function $H:[0,1]\times [0,1]\rightarrow \mathbb D_{\epsilon_1}\times S_0$ with the following properties:

\smallskip

\noindent a) $H(t, 0)=\sigma(t)$, for all $t\in[0,1]$;

\smallskip

\noindent b) $H(t, 1)\in \mathbb D_{\epsilon_1}\times\{y_0\}$, for all $t\in[0,1]$;

\smallskip

\noindent c) $H(0, u) = (x_0, y_0)$, for all $u\in[0,1]$;

\smallskip

\noindent d) if $H(t,u)\in D$, then $t=1$;

\smallskip

\noindent e) $H(1,u)\in C$, for all $u\in[0,1]$;

\smallskip

\noindent f) the path $\eta$ given by $\eta(u):=H(1,u)$ is a homeomorphism onto its image.
\medskip

Thus, $H$ is a homotopy from $\sigma$ to the path $\gamma$ given by $\gamma(t):=H(t, 1)\in \mathbb D_{\epsilon_1}\times\{y_0\}$, such that $(x_0, y_0)$ is ``fixed'' and the point $(x_0, y_1)=H(1, 0)$ ``swings up to the point'' $H(1,1)$ by ``sliding along'' $C$, while the remainder of $\sigma$ does not hit $D$ as it ``swings up'' to $\gamma$.

\end{lem}
\begin{proof} The proper stratified submersion $\mathbb D_{\epsilon_1}\times {\stackrel{{}\hskip -.06in \circ}{\mathbb D^*_{\epsilon_2}}}\arrow{y} {\stackrel{{}\hskip -.06in \circ}{\mathbb D^*_{\epsilon_2}}}$ is a locally trivial fibration, where the local trivialization respects the strata. The restriction of this fibration $\mathbb D_{\epsilon_1}\times S_0\arrow{y} S_0$ is a locally trivial fibration over a contractible space and, hence, is equivalent to the trivial fibration. 

Therefore, there exists a homeomorphism 
$$\Psi: \big(\mathbb D_{\epsilon_1}\times S_0, (\mathbb D_{\epsilon_1}\times S_0)\cap C\big)\rightarrow \big(\mathbb D_{\epsilon_1}\times \{y_0\}, (\mathbb D_{\epsilon_1}\times \{y_0\})\cap C\big)\times[0,1]$$
such that the projection of $\Psi(x, \sigma_0(t))$ onto the $[0,1]$ factor is simply $t$, and such that $\Psi(x, y_0)= ((x, y_0), 0)$. It follows that there is a path $\alpha:[0,1]\rightarrow \mathbb D_{\epsilon_1}$ such that $\Psi(\sigma(t)) = ((\alpha(t), y_0), t)$, for all $t\in[0,1]$. Define $H:[0,1]\times [0,1]\rightarrow \mathbb D_{\epsilon_1}\times S_0$ by
$$
H(t, u):= \Psi^{-1}\big((\alpha(t), y_0), (1-u)t\big).
$$
All of the given properties of $H$ are now trivial to verify.
\end{proof}

\bigskip

\begin{rem}\label{rem:triangle} By Property c) of \lemref{lem:swing}, the map $H$ yields a corresponding map $H^T$ whose domain is a triangle instead of a square. One pictures the image of $H$, or of $H^T$, as a ``gluing in'' of this triangle into $\mathbb D_{\epsilon_1}\times S_0$ in such a way that one edge of the triangle is glued diffeomorphically to $S$, and another edge is glued diffeomorphically onto the image of $\eta$. The third edge of the triangle is glued onto the image of $\gamma$, but not necessarily in a one-to-one fashion.
\end{rem}

\section{The Main Theorem}\label{sec:main} 

In this section, we will prove the Main Theorem.  We must first describe the machinery that goes into this part of the proof.

\bigskip

We will first prove a $1$-dimensional version of the Main Theorem. Assume for now that $\dm_{\mathbf 0}\Sigma f=1$. As the value of $\lambda^1_{f, z_0}(\mathbf 0)$ is minimal for generic $z_0$, we lose no generality if we assume that our linear form $z_0$ is chosen more generically than simply being prepolar. We choose $z_0$ so generically that, in addition to being prepolar, the discriminant, $D$, of the map $(z_0, f)$ and the corresponding Cerf diagram, $C$, have the usual properties -- as given, for instance, in \cite{leperron}, \cite{tibarbouq}, and \cite{vannierthesis}. We will describe the needed properties below.

\bigskip

Let $\widetilde \Psi:=(z_0, f):(\U, \mathbf 0)\rightarrow (\C^2, \mathbf 0)$. We use the coordinates $(u, v)$ on $\C^2$. The critical locus $\Sigma \widetilde \Psi$ of $\widetilde \Psi$ is the union of $\Sigma f$ and $\pc$. The discriminant $D:=\widetilde \Psi(\Sigma \widetilde \Psi)$ consists of the $u$-axis together with the Cerf diagram $C:=\overline{D-V(v)}$. We assume that $z_0$ is generic enough so that the polar curve is reduced and that, in a neighborhood of the origin, $\widetilde \Psi_{|_\pc}$ is one-to-one.

We choose real numbers $\epsilon$, $\delta$, and $\omega$ so that $0< \omega\ll \delta\ll\epsilon\ll 1$. Let $B_\epsilon\subseteq\C^n$ be a closed ball, centered at the origin, of radius $\epsilon$. Let ${\stackrel{\circ}{\D}}_\delta$ and ${\stackrel{\circ}{\D}}_\omega$ be open disks in $\C$, centered at $0$, of radii $\delta$ and $\omega$, respectively.

One considers the map from $({\stackrel{\circ}{\D}}_\delta\times B_\epsilon)\cap f^{-1}({\stackrel{\circ}{\D}}_\omega)$ onto ${\stackrel{\circ}{\D}}_\delta\times{\stackrel{\circ}{\D}}_\omega$ given by the restriction of $\widetilde  \Psi$; we let $ \Psi$ denote this restriction. As $B_\epsilon$ is a closed ball, the map $ \Psi$ is certainly proper, but the domain has an interior stratum, and a stratum coming from the boundary of $B_\epsilon$. However, for generic $z_0$, all of the stratified critical points lie on $\pc\cup\Sigma f$, i.e., above $D$.  

We continue to write simply $D$ and $C$, in place of $D\cap({\stackrel{\circ}{\D}}_\delta\times{\stackrel{\circ}{\D}}_\omega)$ and $C\cap({\stackrel{\circ}{\D}}_\delta\times{\stackrel{\circ}{\D}}_\omega)$. As $ \Psi$ is a proper stratified submersion above ${\stackrel{\circ}{\D}}_\delta\times{\stackrel{\circ}{\D}}_\omega- D$, and as $ \Psi_{|_\pc}$ is one-to-one, many homotopy arguments in $({\stackrel{\circ}{\D}}_\delta\times B_\epsilon)\cap f^{-1}({\stackrel{\circ}{\D}}_\omega)$ can be obtained from lifting constructions in ${\stackrel{\circ}{\D}}_\delta\times{\stackrel{\circ}{\D}}_\omega$. This is the point of considering the discriminant and Cerf diagram. 

Let $v_0\in{\stackrel{\circ}{\D}}_\omega-\{\mathbf 0\}$. By construction, up to diffeomorphism, $ \Psi^{-1}({\stackrel{\circ}{\D}}_\delta\times\{v_0\})$ is $\mf$ and $ \Psi^{-1}((0, v_0))$ is $\mfoo$. In fact, for all $u_0$, where $|u_0|\ll |v_0|$, $ \Psi^{-1}((u_0, v_0))$ is diffeomorphic to $\mfoo$; we fix such a non-zero $u_0$, and let $\mathbf a:=(u_0, v_0)$.

\smallskip

We wish to pick a distinguished basis for the vanishing cycles of $f_0$ at the origin, as in I.1 of \cite{agv} (see, also, \cite{dimcasing}). We do this by selecting paths in $\{u_0\}\times {\stackrel{\circ}{\D}}_\omega$ which originate at $\mathbf a$. We must be slightly careful in how we do this. 

First, fix a path $p_0$ from $\mathbf a$ to $(u_0, 0)$. Select paths  $q_1, \dots, q_{\gamma^1}$ from $\mathbf a$ to each of the points in $(\{u_0\}\times {\stackrel{\circ}{\D}}_\omega)\cap C=:\{y_1, \dots, y_{\gamma^1}\}$. The paths $p_0, q_1, \dots, q_{\gamma^1}$ should not intersect each other and should intersect the set $\{(u_0, 0), y_1, \dots, y_{\gamma^1}\}$  only at the endpoints of the paths. Moreover, when at the point $\mathbf a$, the paths $p_0, q_1, \dots, q_{\gamma^1}$ should be in clockwise order. Let $r_0$ be a clockwise loop very close to $p_{0}$, from $\mathbf a$ around $(u_0, 0)$.

As we are not assuming that $f$ had an isolated line singularity, we must perturb $f_{|_{V(z_0-u_0)}}$ slightly to have $(u_0, 0)$ split into $\lambda^1$ points, $x_1, \dots, x_{\lambda^1}$ inside the loop $r_0$; each of these points corresponds to an $A_1$ singularity in the domain. We select paths $p_1, \dots, p_{\lambda^1}$ from $\mathbf a$ to each of the points $x_1, \dots, x_{\lambda^1}$, and paths $q_1, \dots, q_{\gamma^1}$ from $\mathbf a$ to each of the points in $(\{u_0\}\times {\stackrel{\circ}{\D}}_\omega)\cap C=:\{y_1, \dots, y_{\gamma^1}\}$. We may do this in such a way that the paths $p_1, \dots, p_{\lambda^1}, q_1, \dots, q_{\gamma^1}$ are in clockwise order. 

\smallskip

The lifts of these paths via the perturbed $f_{|_{V(z_0 - u_0)}}$ yield representatives of elements of $H_{n+1}(B_\epsilon, \mfoo)$, whose boundaries in $\widetilde H_n(\mfoo)$ form a distinguished basis  $\Delta_1^\prime, \dots, \Delta_{\lambda^1}^\prime, \Delta_1, \dots, \Delta_{\gamma^1}$.

\smallskip

By using the swing (\lemref{lem:swing}), the paths $q_1, \dots, q_{\gamma^1}$ are taken to new paths $\hat q_1, \dots, \hat q_{\gamma^1}$ in ${\stackrel{\circ}{\D}}_\delta\times\{v_0\}$. Each $\hat q_i$ path represents a relative homology class in $H_n(\mf, \mfoo)$ whose boundary in $\widetilde H_{n-1}(\mfoo)$ is precisely $\Delta_i$.  \thmref{thm:swing} follows from this. 

\bigskip

We can now prove:

\bigskip

\begin{thm}\label{thm:mainone} Suppose that $\dm_{\mathbf 0}\Sigma f=1$ and $\dm_{\mathbf 0}\Sigma f_0=0$. Then, the following are equivalent:

\vskip .1in

\noindent a) $f_{\mathbf q}$ is a simple $\mu$-constant family, i.e., $f$ has a smooth critical locus which defines a family of isolated singularities with constant Milnor number $\mu_{f_0}$;

\vskip .1in

\noindent b) $\rank\widetilde H_{n-1}(\mf)=\lambda^1_{f, z_0}(\mathbf 0)$;

\vskip .1in

\noindent c) there exists a prime $p$ such that $\dm \widetilde H_{n-1}(\mf;\ \Z/p\Z)=\lambda^1_{f, z_0}(\mathbf 0)$.

\vskip .1in

Thus,if $f$ is not Milnor equisingular, $\rank\widetilde H_{n-1}(\mf)<\lambda^1_{f, z_0}(\mathbf 0)$, and so $\rank\widetilde H_{n}(\mf)<\lambda^0_{f, z_0}(\mathbf 0)$, and these inequalities hold with $\Z/p\Z$ coefficients (here, $p$ is prime).\end{thm}

\begin{proof} That a) implies b) and c) is well-known; it follows at once from \thmref{thm:leattach}. Assume then that $f_{\mathbf q}$ is not a simple $\mu$-constant family. We will prove that $\rank\widetilde H_{n-1}(\mf)<\lambda^1_{f, z_0}(\mathbf 0)$, and then indicate why the same proof applies with $\mathbb Z/p\mathbb Z$ coefficients.

By \thmref{thm:milnorequi}, $\pc\neq\emptyset$, and so $C\neq\emptyset$. We want to construct just one new path in  $\{u_0\}\times {\stackrel{\circ}{\D}}_\omega$, one which originates at $\mathbf a$, ends at a point of $C$, and misses all of the other points of $D$; we want this path to swing up to a path in  ${\stackrel{\circ}{\D}}_\delta\times\{v_0\}$, and represent a relative homology class in $H_n(\mf, \mfoo)$ whose boundary is not in the span of $\Delta_1, \dots, \Delta_{\gamma^1}$.

By the connectivity of the vanishing cycle intersection diagram (\cite{gabrielov}, \cite{lazzeri}), one of the $\Delta^\prime_j$ must have a non-zero intersection pairing with one of the $\Delta_i$, i.e., there exist $i_0$ and $j_0$ such that $\langle \Delta_{i_0}, \Delta^\prime_{j_0}\rangle\neq 0$. 

By fixing the path $p_{j_0}$ and all the $q_i$ paths, but reselecting the other $p_j$, for $j\neq j_0$, we may assume that $j_0=1$, i.e., that $\langle \Delta_{i_0}, \Delta^\prime_{1}\rangle\neq 0$. 

We follow now Chapter 3.3 of \cite{dimcasing}. Associated to each path $p_{j}$, $1\leqslant j\leqslant\lambda^1$, is a (partial) monodromy automorphism $T^\prime_{j}:\widetilde H_{n-1}(\mfoo)\rightarrow \widetilde H_{n-1}(\mfoo)$, induced by taking a clockwise loop $r_{j}$ very close to $p_{j}$, from $\mathbf a$ around $x_{j}$. Let $T^\prime:= T^\prime_1\dots T^\prime_{\lambda^1}$, where composition is written in the order of \cite{dimcasing}. We claim that $T^\prime(\Delta_{i_0})$ is in the image of $\delta: H_n(\mf, \mfoo)\rightarrow \widetilde H_{n-1}(\mfoo)$, but is not in $\operatorname{Span}\{\Delta_1, \dots, \Delta_{\gamma^1}\}$.

The composition $r$ of the loops $r_1, \dots, r_{\lambda^1}$ is homotopy-equivalent, in $\{u_0\}\times {\stackrel{\circ}{\D}}_\omega-\big\{\{x_1, \dots, x_{\lambda^1}\}\cup C\big\}$, to the loop $r_0$ (from our discussion before the theorem).   By combining (concatenating) the loop $r_{0}$ and the path $q_{i_0}$, we obtain a path in $\{u_0\}\times {\stackrel{\circ}{\D}}_\omega$ which is homotopy-equivalent to a simple path which swings up to a corresponding path in  ${\stackrel{\circ}{\D}}_\delta\times\{v_0\}$. Thus, $T^\prime(\Delta_{i_0})$ is in the image of $\delta$.

Now, by the Corollaries to the Picard-Lefschetz Theorem in \cite{agv}, p. 26, or as in \cite{dimcasing}, Formula 3.11,
$$
T^\prime(\Delta_{i_0}) = \Delta_{i_0}- (-1)^{\frac{n(n-1)}{2}}\langle \Delta_{i_0}, \Delta^\prime_{1} \rangle \Delta^\prime_{1}+\beta_2 \Delta^\prime_{2}+\dots+\beta_{\lambda^1} \Delta^\prime_{\lambda^1},
$$
for some integers $\beta_2, \dots, \beta_{\lambda^1}$.
As the $\Delta_1^\prime, \dots, \Delta_{\lambda^1}^\prime, \Delta_1, \dots, \Delta_{\gamma^1}$ form a basis, and as $\langle \Delta_{i_0}, \Delta^\prime_{1} \rangle\neq 0$, $T^\prime(\Delta_{i_0})$ is not in $\operatorname{Span}\{\Delta_1, \dots, \Delta_{\gamma^1}\}$.

This finishes the proof over the integers. Over $\Z/p\Z$, the proof is identical, since the intersection diagram is also connected modulo $p$; see \cite{gabrielov}.
\end{proof}

\bigskip

\begin{rem}\label{rem:careful} One must be careful in the proof above; it is tempting to try to use simply $T^\prime_1(\Delta_{i_0})$ in place of $T^\prime(\Delta_{i_0})$. The problem with this is that $T^\prime_1(\Delta_{i_0})$ is not represented by a path in $\{u_0\}\times {\stackrel{\circ}{\D}}_\omega - \{(u_0, 0)\}$ and, thus, there is no guaranteed swing isotopy to a corresponding path in  ${\stackrel{\circ}{\D}}_\delta\times\{v_0\}$.

In fact, we could have avoided the explicit construction of   $T^\prime(\Delta_{i_0})$ completely, though we find the construction intuitive and geometrically interesting. By naturality (of the monodromy automorphism on the vanishing cycle functor), the map $\delta: H_n(\mf, \mfoo)\rightarrow \widetilde H_{n-1}(\mfoo)$ commutes with the respective monodromy actions. Thus, the image of $\delta$, $\operatorname{im}\delta$, is invariant under the monodromy action. Now, the swing and the construction of the distinguished basis for $\widetilde H_{n-1}(\mfoo)$ tell us that we can write $\widetilde H_{n-1}(\mfoo)$ as a direct sum $A\oplus B$, where $A$ and $B$ are generated by distinguished basis elements, $\operatorname{rank} A = \gamma^1$, and $A\subseteq \operatorname{im}\delta$. However, the connectivity of the intersection matrix for $f_0$ implies that the only monodromy-invariant submodules of $\widetilde H_{n-1}(\mfoo)$, which are generated by distinguished basis elements, are the zero-module and all of $\widetilde H_{n-1}(\mfoo)$ (see \cite{agv}, Theorem 3.5). Thus, if $\gamma^1\neq 0$ (i.e., if we do not have a simple $\mu$-constant family), then the image of $\delta$ has to properly contain $A$. \thmref{thm:mainone} follows.
\end{rem}

\vskip .4in

We can now prove our Main Theorem. We return to the general case where $s:=\dm_{\mathbf 0}\Sigma f$ is arbitrary. Fix a set of coordinates $(z_0, \dots, z_n)$, and consider the corresponding family $f_{\mathbf q}$.

\bigskip

\begin{thm}{\rm \bf(Main Theorem)}\label{thm:main}.  Suppose that $\dm_{\mathbf 0}\Sigma(f_{\mathbf 0}) = 0$.

Then, $\rank\widetilde H_{n-s}(\mf)=\lambda^s_{f, \mathbf z}(\mathbf 0)$ if and only if $f_{\mathbf q}$ is a simple $\mu$-constant family.
\end{thm}

\begin{proof} If $f_{\mathbf q}$ is a simple $\mu$-constant family, then it is well-known that $\rank\widetilde H_{n-s}(\mf)=\lambda^s_{f, \mathbf z}(\mathbf 0)$; this follows from an inductive application of \cite{leattach}, using that $\Gamma^s_{f, \mathbf z}=0$ (as we saw in \thmref{thm:milnorequi}).

Suppose that $\tilde b_{n-s}:=\rank\widetilde H_{n-s}(\mf)=\lambda^s_{f, \mathbf z}(\mathbf 0)$. As we saw in the Introduction, $\tilde b_{n-s}\leq \lambda^s_f(\mathbf 0)$. Thus, we must have that $\lambda^s_{f, \mathbf z}(\mathbf 0)= \lambda^s_f(\mathbf 0)$. Let $\hat \mathbf z$ be a generic choice of coordinates at the origin, consider the codimension $s-1$ linear slice $N:=V(\hat z_0, \dots, \hat z_{s-2})$ through the origin. 

Then, $f_{|_N}$ has a $1$-dimensional critical locus and, by iterating \thmref{thm:leattach}, $\widetilde H_{n-s}(\mf)\cong \widetilde H_{(n-s+1)-1}(F_{f_{|_N}})$. Now, by Proposition 1.21 of \cite{lecycles}, $\lambda^s_{f, \hat\mathbf z}(\mathbf 0) = \lambda^1_{f_{|_N}, \hat z_{s-1}}(\mathbf 0)$.  \thmref{thm:main} implies that $f_{|_N}$ is Milnor equisingular; in particular, the polar curve, $\Gamma^1_{f_{|_N}, \hat z_{s-1}}$ is zero (or, as a set, is empty). By Proposition 1.21 of \cite{lecycles}, this implies that $\Gamma^s_{f, \hat\mathbf z}=0$. Now, by d) of \thmref{thm:milnorequi}, $f$ is Milnor equisingular. Therefore, $\Sigma f$ is smooth at the origin and, since $\lambda^s_{f, \mathbf z}(\mathbf 0)= \lambda^s_f(\mathbf 0)$, $V(z_0, \dots, z_{s-1})$ must transversely intersect $\Sigma f$ at the origin. The desired conclusion now follows from e) of \thmref{thm:milnorequi}.
\end{proof}

\bigskip

We wish to see that \thmref{thm:main} puts restrictions on the types of perverse sheaves that one may obtain as vanishing cycles $\phi_f[-1]\mathbb Z^\bullet_{\U}[n+1]$ of the shifted constant sheaf on affine space. Below, we refer to the constant sheaf on $\nu$ of rank ${\stackrel{\circ}{\mu}}_\nu$, shifted by some integer $j$ and extended by zero to all of $V(f)$; we write $(\mathbb Z^{{\stackrel{\circ}{\mu}}_\nu})^\bullet_\nu[j]$ for this sheaf (note that we omit the reference to the extension by zero in the notation). The isomorphisms and direct sums that we write below are in the Abelian category of perverse sheaves.

\bigskip

\begin{cor}\label{cor:pervstruct} Suppose that the critical locus of $f$ is  $s$-dimensional, where $s\geqslant 1$ is arbitrary, and that every $s$-dimensional component, $\nu$, of $\Sigma f$ is smooth.

Then,  $\bigoplus_\nu(\mathbb Z^{{\stackrel{\circ}{\mu}}_\nu})^\bullet_\nu[s]$  is a direct summand of $\phi_f[-1]\mathbb Z^\bullet_{\U}[n+1]$ if and only if $f$ is Milnor equisingular. Moreover, when these equivalent conditions hold, $\Sigma f$ is smooth and $\phi_f[-1]\mathbb Z^\bullet_{\U}[n+1]\cong (\mathbb Z^{{\stackrel{\circ}{\mu}}_{\Sigma f}})^\bullet_{{}_{\Sigma f}}[s]$.
\end{cor}

\begin{proof}
If $f$ is Milnor equisingular, then, by b) of \thmref{thm:milnorequi}, $V(f)$ has an $a_f$ stratification consisting of two strata: $V(f)-\Sigma f$ and $\Sigma f$. As $\phi_f[-1]\mathbb Z^\bullet_{\U}[n+1]$ is constructible with respect to any $a_f$ stratification, it follows that $\phi_f[-1]\mathbb Z^\bullet_{\U}[n+1]\cong (\mathbb Z^{{\stackrel{\circ}{\mu}}_{\Sigma f}})^\bullet_{{}_{\Sigma f}}[s]$.

\smallskip

If $\bigoplus_\nu(\mathbb Z^{{\stackrel{\circ}{\mu}}_\nu})^\bullet_\nu[s]$  is a direct summand of $\phi_f[-1]\mathbb Z^\bullet_{\U}[n+1]$, then $\rank\widetilde H_{n-s}(\mf)$ is at least $\sum_\nu {\stackrel{\circ}{\mu}}_\nu$, which equals $\lambda^s_f(\mathbf 0)$, as each $\nu$ is smooth. Now, \thmref{thm:main} tells us that $f$ must be Milnor equisingular.
\end{proof}

\section{Comments, Questions, and Counterexamples}\label{sec:cqc}

 One might hope that a stronger result than \thmref{thm:main}, or \thmref{thm:mainone}, is true.
 
\bigskip

For instance, given that \thmref{thm:mainone} and \thmref{thm:siersmabound} are true, it is natural to ask the following:

\begin{ques} If we are not in the trivial case, is the rank of $\widetilde H_{n-1}(\mf)$ strictly less than $\sum_\nu {\stackrel{\circ}{\mu}}_\nu$?
\end{ques}

\smallskip

The answer to the above question is ``no''. One can find examples of this in the literature, but perhaps the easiest is the following:

\smallskip

\begin{exm} Let $f:= (y^2-x^3)^2+w^2$. Then, $\Sigma f$ has a single component $\nu:=V(w, y^2-x^3)$, and one easily checks that ${\stackrel{\circ}{\mu}}_\nu = 1$. However, as $f$ is the suspension of $(y^2-x^3)^2$, the Sebastiani-Thom Theorem (here, we need the version proved by Oka in \cite{okasebthom}) implies
$$
\widetilde H_1(\mf)\cong\widetilde H_0(F_{(y^2-x^3)^2})\cong\Z.
$$
Moreover, by suspending $f$ again, one may produce an example in which $f$ itself has a single irreducible component at the origin.

It is not difficult to show that, for this example, ${\mathbb Z}^\bullet_{{}_{\Sigma f}}[1]$  is a direct summand of $\phi_f[-1]\mathbb Z^\bullet_{\U}[n+1]$. Thus, this example shows that the assumption on the smoothness of the $s$-dimensional components of $\Sigma f$ in \corref{cor:pervstruct} is necessary. This is especially interesting since  $\Sigma f$ is homeomorphic to a complex line, and perverse sheaves are topological devices. This shows that the structure of the vanishing cycles in the perverse category ``remembers'' the hypersurface that surrounded $\Sigma f$.
\end{exm}

\bigskip

 Now, let $\alpha$ be the number of irreducible components of $\Sigma f$.
 
 \smallskip
 
 \begin{ques}\label{ques:betterbound} If we are not in the trivial case, is the rank of $\widetilde H_{n-1}(\mf)$ strictly less than $\lambda^1-\alpha$?
\end{ques}

\smallskip

Again, there are many examples in the literature which demonstrate that the answer to this question is ``no''. One simple example is:

\begin{exm} 
The function $f=x^2y^2+w^2$ has a critical locus consisting of two lines, $\lambda^1=2$, but -- using the Sebastiani-Thom Theorem again -- we find that $\widetilde H_1(\mf)\cong\Z$.
\end{exm}

\bigskip

However, a result such as that asked about in \quesref{ques:betterbound}, but where $\alpha$ is replaced by a quantity involving the number of components of $\pc$, or numbers of various types of components in the Cerf diagram, seems more likely. Moreover, if we put more conditions on the intersection diagram for the vanishing cycles of $f_0$, we could certainly obtain sharper bounds than we do in the Main Theorem.  Or, if we know more topological data, such as the vertical monodromies, as in \cite{siersmavarlad}, we could obtain better bounds. However, other than \thmref{thm:mainone}  and  \thmref{thm:main}, we know of no nice, effectively calculable, bound which holds in all cases.

\bigskip

Finally, \corref{cor:pervstruct} leads us to ask: 

\smallskip

\begin{ques}\label{ques:pervstruct} Which perverse sheaves can be obtained as the vanishing cycles of the constant sheaf on affine space?
\end{ques}

\bigskip

Unlike our previous questions, we do not know the answer to \quesref{ques:pervstruct}.

\bibliographystyle{plain}
\bibliography{Masseybib}

\begin{thebibliography}{10}

\bibitem{agv}
{Arnold, V. I., Gusein-Zade, S. M., Varchenko, A. N.}
\newblock {\em {Singularities of Differentiable Maps II, Monodromy and
  Asymptotics of Integrals}}.
\newblock {Birkh\"auser}, 1988.

\bibitem{caubelthesis}
{Caubel, C.}
\newblock {\em {Sur la topologie d'une famille de pinceaux de germes
  d'hypersurfaces complexes}}.
\newblock PhD thesis, Universit\'e Toulouse III, 1998.

\bibitem{dejong}
{de Jong, Th.}
\newblock {Some classes of line singularities}.
\newblock {\em Math. Zeitschrift}, 198:493--517, 1998.

\bibitem{dimcasing}
{Dimca, A.}
\newblock {\em {Singularities and Topology of Hypersurfaces}}.
\newblock Universitext. Springer-Verlag, 1992.

\bibitem{fulton}
{Fulton, W.}
\newblock {\em {Intersection Theory}}.
\newblock Ergeb. Math. Springer-Verlag, 1984.

\bibitem{gabrielov}
{Gabrielov, A. M.}
\newblock {Bifurcations, Dynkin Diagrams, and Modality of Isolated
  Singularities}.
\newblock {\em Funk. Anal. Pril.}, 8 (2):7--12, 1974.

\bibitem{katomatsu}
{Kato, M. and Matsumoto, Y.}
\newblock {On the connectivity of the Milnor fibre of a holomorphic function at
  a critical point}.
\newblock {\em Proc. of 1973 Tokyo manifolds conf.}, pages 131--136, 1973.

\bibitem{lazzeri}
{Lazzeri, F.}
\newblock {\em {Some Remarks on the Picard-Lefschetz Monodromy}}.
\newblock {Quelques journ\'ees singuli\`eres}. Centre de Math. de l'Ecole
  Polytechnique, Paris, 1974.

\bibitem{leacampo}
{L\^e, D. T. }.
\newblock {Une application d'un th\'eor\`eme d'A'Campo a l'equisingularit\'e}.
\newblock {\em Indag. Math}, 35:403--409, 1973.

\bibitem{leattach}
{L\^e, D. T.}
\newblock {Calcul du Nombre de Cycles Evanouissants d'une Hypersurface
  Complexe}.
\newblock {\em Ann. Inst. Fourier, Grenoble}, 23:261--270, 1973.

\bibitem{leperron}
{L\^e, D. T. and Perron, B.}
\newblock {Sur la fibre de Milnor d'une singularit\'e isol\'ee en dimension
  complexe trois}.
\newblock {\em C. R. Acad. Sci. Pairs S\'er. A}, 289:115--118, 1979.

\bibitem{leramanujam}
{L\^e, D. T. and Ramanujam, C. P.}
\newblock {The Invariance of Milnor's Number implies the Invariance of the
  Topological Type}.
\newblock {\em Amer. Journ. Math.}, 98:67--78, 1976.

\bibitem{lesaito}
{L\^e, D. T. and Saito, K.}
\newblock {La constance du nombre de Milnor donne des bonnes stratifications}.
\newblock {\em C.R. Acad. Sci.}, 277:793--795, 1973.

\bibitem{lecycles}
{Massey, D.}
\newblock {\em {L\^e Cycles and Hypersurface Singularities}}, volume 1615 of
  {\em Lecture Notes in Math.}
\newblock Springer-Verlag, 1995.

\bibitem{okasebthom}
{Oka, M.}
\newblock {On the homotopy type of hypersurfaces defined by weighted
  homogeneous polynomials}.
\newblock {\em Topology}, 12:19--32, 1973.

\bibitem{siersmaisoline}
{Siersma, D.}
\newblock {Isolated line singularities}.
\newblock {\em Proc. Symp. Pure Math.}, 35, part 2, Arcata Singularities
  Conf.:485--496, 1983.

\bibitem{siersmavarlad}
{Siersma, D.}
\newblock {Variation mappings on singularities with a $1$-dimensional critical
  locus}.
\newblock {\em Topology}, 30:445--469, 1991.

\bibitem{teissiercargese}
{Teissier, B.}
\newblock {Cycles \'evanescents, sections planes et conditions de Whitney}.
\newblock {\em Ast\'erisque}, 7-8:285--362, 1973.

\bibitem{tibarbouq}
{Tib\u ar, M.}
\newblock {Bouquet Decomposition of the Milnor Fiber}.
\newblock {\em Topology}, 35:227--241, 1996.

\bibitem{vannierthesis}
{Vannier, J. P.}
\newblock {\em {Familles \`a param\`etre de fonctions holomorphes \`a ensemble
  singulier de dimension z\'ero ou un}}.
\newblock PhD thesis, Dijon, 1987.

\end{thebibliography}
\end{document}